\magnification 1200
\input amstex
\documentstyle{amsppt}
\pageheight{19cm}
\NoBlackBoxes
\define\al{\alpha}
\define\be{\beta}
\define\ga{\gamma}
\define\de{\delta}
\define\la{\lambda}

\define\Ga{\Gamma}

\define\CC{{\Bbb C}}
\define\ZZ{{\Bbb Z}}

\define\const{\text{const.\,}}
\define\rep{representation}
%
\topmatter
\title
On $q$-analogues of the Fourier and Hankel transforms
\endtitle
\author
TOM H. KOORNWINDER and
REN\'E F. SWARTTOUW
\endauthor
\address
CWI, Amsterdam
\endaddress
\curraddr
University of Amsterdam, Korteweg-de Vries Institute,
P.O.Box 94248,\quad 1090 GE Amsterdam, The Netherlands
\endcurraddr
\email
T.H.Koornwinder\@uva.nl
\endemail
\address
Faculty of Technical Mathematics and Informatics,
Delft University of Technology
\endaddress
\curraddr
VU University Amsterdam,
Faculty of Exact Sciences,
Department of Mathematics,
De Boelelaan 1081,
1081 HV Amsterdam,
The Netherlands
\endcurraddr
\email
rene\@few.vu.nl
\endemail
\thanks
This paper appeared in Trans.\ Amer.\ Math.\ Soc. 333 (1992), 445--461.
In the present version the statement and proof of Proposition 4.2
and the proof of Proposition A.1 have been corrected.
\endthanks
\abstract
For the third $q$-Bessel function
(first introduced by F.~H. Jackson, later rediscovered by W.~Hahn in
a special case and by H. Exton)
we derive Hansen-Lommel
type orthogonality relations, which, by a symmetry, turn out to be
equivalent to orthogonality relations which are $q$-analogues of the
Hankel integral transform pair. These results are implicit, in the context
of quantum groups, in a paper by Vaksman and Korogodski{\u\i}.
As a specialization we get $q$-cosines
and $q$-sines which admit $q$-analogues of the Fourier-cosine and Fourier-sine
transforms. We get also a formula which is both an analogue of Graf's
addition formula and of the Weber-Schafheitlin discontinuous integral.
This is a corrected version of a paper which appeared in
Trans.\ Amer.\ Math.\ Soc. 333 (1992), 445--461.
\endabstract
\endtopmatter
\document

\heading 1. Introduction\endheading
Several possible $q$-analogues of the {\sl Bessel function}
$$J_\al(x):=\sum_{k=0}^\infty
{(-1)^k\,(x/2)^{\al+2k}\over k!\,\Ga(\al+k+1)}
\tag 1.1$$
have been considered in the literature.
Best known are two related $q$-Bessel functions denoted
$J_\al^{(1)}(x;q)$ and $J_\al^{(2)}(x;q)$ by
Ismail \cite{10},
but first introduced by Jackson in a series of papers during the years
1903--1905 (see the references in \cite{10})
and also studied by Hahn \cite{7}.
A third $q$-Bessel function was introduced by Hahn \cite{8}
(in a special case; we thank G.~Gasper for this reference) and by
Exton \cite{3}, \cite{4, (5.3.1.11)} (in full). In Exton's
notation:
$$C_\al(q;x):=
{(1-q)^\al\over (q;q)_\infty}\,
\sum_{k=0}^\infty{q^{k(k-1)/2}\,(q^{\al+k+1};q)_\infty\,(-x(1-q)^2)^k
\over(q;q)_k}\,,
\tag 1.2$$
where the $q$-shifted factorials are defined by
$$(a;q)_k:=\prod_{j=0}^{k-1}(1-aq^j),\qquad k=1,2,\dots;\qquad(a;q)_0:=1;
\tag 1.3$$
$$(a;q)_\infty:=\lim_{k\to\infty}(a;q)_k,\qquad|q|<1.$$
Hahn \cite{8} considered the case $\al=0$ of (1.2).
They obtained these functions as the solutions of a special basic
Sturm-Liouville equation, by which they could also derive the following
$q$-analogue of the Fourier-Bessel orthogonality relations:
$$\int_0^1 x^\al\,C_\al(q;\mu_i qx)\,C_\al(q;\mu_j qx)\,d_q x=0,\qquad
i\ne j,\quad\al>-1,$$
where $\mu_1,\mu_2,\dotsc$ are the roots of the equation
$C_\al(q;\mu)=0$,
and where the $q$-integral is defined by
$$\int_0^1 f(x)\,d_q x:=(1-q)\,\sum_{k=0}^\infty f(q^k)\,q^k.$$
By specialization to $\al=\pm1/2$, Exton obtained similar orthogonalities
for $q$-analogues of the sines and cosines.
So some of the harmonic analysis involving Bessel functions, sines and
cosines has been extended to the $q$-case. However, $q$-analogues of the
Fourier-cosine, Fourier-sine and Hankel transforms were missing until now
(except for a $q$-Laplace transform with inversion formula,
cf\. Hahn \cite{7, \S9} and Feinsilver \cite{5}).

Exton \cite{4, (5.3.3.1)}
also generalized the generating function
$$e^{z(t-t^{-1})/2}=\sum_{k=-\infty}^\infty t^k \, J_k(z)
\tag 1.4$$
(cf\. Watson \cite{15, \S2.1 (1)})
to the case of his $q$-Bessel functions (1.2):
$$e_q((1-q)t)\,E_q(-(1-q)t^{-1}x)=\sum_{n=-\infty}^\infty t^n\,C_n(q;x),
\tag 1.5$$
where
$$e_q(z):=\sum_{k=0}^\infty {z^k\over(q;q)_k}={1\over(z;q)_\infty}
\, ,\qquad|z|<1,\tag 1.6$$
and
$$E_q(z):=\sum_{k=0}^\infty{q^{k(k-1)/2}\,z^k\over(q;q)_k}=(-z;q)_\infty
\tag 1.7$$
are $q$-analogues of the exponential function
(cf\. Gasper and Rahman \cite{6, (1.3.15), (1.3.16)}).

Recently, Vaksman and Korogodski{\u\i} \cite{14} gave an
interpretation
of the $q$-Bessel functions (1.2) as matrix elements of irreducible
\rep s
of the quantum group of plane motions. Their paper, which does not contain
proofs, implicitly contains some new orthogonality relations for the functions
(1.2). In particular, the unitariness of the \rep s implies a
$q$-analogue of the Hansen-Lommel orthogonality relations
$$\de_{nm}=\sum_{k=-\infty}^\infty J_{k+n}(x)\,J_{k+m}(x),\qquad n,m\in\ZZ
\tag 1.8$$
(cf\. \cite{15, \S2.5 (3),(4)}).
Furthermore, the Schur type orthogonality relations for matrix elements of
irreducible unitary \rep s which are square integrable with respect to a
suitable Haar functional, yield a $q$-analogue of Hankel's Fourier-Bessel
integral $$f(x)=\int_0^\infty J_\al(xt)\,
\bigl(\int_0^\infty J_\al(ty)\,f(y)\,y\,dy\bigr)\,t\,dt,
\tag 1.9$$
cf\. \cite{15, \S14.3}.

It is the purpose of the present note to state these two types of orthogonality
relations for the functions (1.2) explicitly, to show that the first
type is immediately implied by the generating function (1.5), and to
rewrite the first type as the second type by use of a simple, but possibly new,
symmetry for the functions (1.2). This will be done in \S2.
In \S3 we will show that the second type of orthogonality is, on the one
hand, a limit case of the orthogonality for the little $q$-Jacobi polynomials,
and, on the other hand, allows the Hankel transform inversion formula as a
limit case for $q\uparrow1$.
In \S4 we will generalize the two orthogonality relations to
two equivalent formulas, which are respectively the $q$-analogues of
Graf's addition formula and the Weber-Schafheitlin discontinuous
integral.
The special cases of the $q$-Fourier-cosine and sine
transforms will be the topic of \S5. No material from \S4 is needed in
this section.
Finally, Appendix A will contain rigorous
proofs of two limit results.

In a recent paper by Rahman \cite{13}, where a $q$-analogue of the
Fourier-Bessel orthogonality for Jackson's $q$-Bessel functions is discussed,
the author states in his concluding remarks that Jackson's $q$-Bessel functions
have probably nicer properties than those of Exton. However, the
results obtained in \cite{14} and in the present paper might suggest
that the Hahn-Exton functions are more suitable for
harmonic analysis, both within and without the context of quantum groups.
Future research will help to clarify the merits of
the various types of $q$-Bessel functions.

In this paper we will not preserve Exton's notation $C_\al(q;x)$ in
(1.2), but state the results in terms of the
$q$-hypergeometric function
$${}_1\phi_1(0;w;q,z):=
\sum_{k=0}^\infty {(-1)^k\,q^{k(k-1)/2}\,z^k\over
(w;q)_k\,(q;q)_k}\,.
\tag 1.10$$
Our motivation is that (i) the ${}_r\phi_s$ notation is fairly well known
nowadays and clarifies the position of these $q$-Bessel functions among other
$q$-hypergeometric functions;
(ii) $q$-analysis should not depend too much on the $q=1$ case, so scaling
factors making the limit transition $q\uparrow1$ easy should not be hidden
in the definitions of $q$-special functions; and
(iii) the classical definition (1.1) of Bessel functions is very natural
in the context of the generating function (1.4) and for analysis
concentrating on $J_0$, but it is less fortunate when the focus is on
another Bessel function of fixed order (cf\. $J_{-1/2}(x)$ versus $\cos x$),
so one should be very careful before fixing the definition and notation of
a $q$-Bessel function. Therefore, the notations for $q$-Bessel functions,
$q$-cosines and $q$-sines in \S\S3 and 5 should be considered as
ad hoc notations, which are only used locally in this paper for clarifying
the analogy of the functions with the $q=1$ case.

While this Introduction until here describes what was known to us in 1992,
it has later been observed by M.~E.~H. Ismail (see for instance p.184 in
\cite{16})
that this third type of $q$-Bessel
function was already observed by Jackson on p.201 of his paper \cite{17}.
His notation and definition is
$$
\text{J}_{[n]}(\lambda):=\sum_{r=0}^\infty (-1)^r p^{r(r-1)}\,
{((1-p)\lambda)^{n+2r}\over(p^2;p^2)_{n+r}(p^2;p^2)_r}\,.
$$
Therefore, what we called the Hahn-Exton $q$-Bessel function is now often called
{\it Jackson's third $q$-Bessel function}.

\heading 2. Symmetry and orthogonality for $q$-Bessel functions\endheading
We will always assume that $0<q<1$.
The general $q$-hypergeometric series is defined by
$${}_r\phi_s\left[{a_1,\dots,a_r\atop b_1,\dots,b_s};q,z\right]:=
\sum_{k=0}^\infty
{(a_1,\dots,a_r;q)_k\,\bigl((-1)^k\,q^{ k(k-1)/2}\bigr)^{s-r+1}\,z^k
\over
(b_1,\dots,b_s;q)_k\,(q;q)_k}\,,
\tag 2.1$$
where the $q$-shifted factorial is defined by (1.3) and
$$(a_1,\dots,a_r;q)_k:=(a_1;q)_k\,(a_2;q)_k\dotsb(a_r;q)_k.$$
The upper and lower parameters in the left-hand side of (2.1) may also be
written on one line as $a_1,\dots,a_r;b_1,\dots,b_s$.
The power series, in the non-terminating case of (2.1), has radius of
convergence $\infty$, 1 or 0 according to whether $r-s<1$, $=1$ or $>1$,
respectively (see \cite{6,  Ch.1} for further details).
Thus the defining formula
(1.2) for the Hahn-Exton $q$-Bessel function can be rewritten as
$$C_\al(q;x):={(1-q)^\al\,(q^{\al+1};q)_\infty\over(q;q)_\infty}\,
{}_1\phi_1(0;q^{\al+1};q, x(1-q)^2)$$
and the $q$-exponential functions (1.6), (1.7)
can be written as
${}_1\phi_0(0;-;q,z)$ and ${}_0\phi_0(-;-;q,\allowbreak-z)$, respectively.

Our object will be the Hahn-Exton $q$-Bessel function written as
the $q$-hypergeometric series (1.10).
It is well-defined for $z,w\in\CC$ with $w$ outside
$\{1,q^{-1},q^{-2},\ldots\}$.
These singularities can be removed by multiplication by $(w;q)_\infty$:
$$(w;q)_\infty\,{}_1\phi_1(0;w;q,z)=
\sum_{k=0}^\infty
{(-1)^k\,q^{k(k-1)/2}\,(q^k w;q)_\infty\,z^k \over (q;q)_k} \, .
\tag 2.2$$

\proclaim{Proposition 2.1}
The series in (2.2) defines an entire analytic function in $z,w$,
which is also symmetric in $z,w$:
$$(w;q)_\infty\,{}_1\phi_1(0;w;q,z)=(z;q)_\infty\,{}_1\phi_1(0;z;q,w).
\tag 2.3$$
Both sides can be majorized by
$$(-|z|;q)_\infty \, (-|w|;q)_\infty .
\tag 2.4$$
\endproclaim
\demo{Proof}
Substitute for $(q^k w;q)_\infty$ in (2.2) the ${}_0\phi_0$ series
given by (1.7):
$$(w;q)_\infty\,{}_1\phi_1(0;w;q,z)=
\sum_{k=0}^\infty\,\sum_{l=0}^\infty
q^{kl}\,
{(-1)^l\,q^{l(l-1)/2}\,w^l \over (q;q)_l}\,
{(-1)^k\,q^{k(k-1)/2}\,z^k \over (q;q)_k} .$$
The summand of the double series can be majorized by
$${q^{l(l-1)/2}\,|w|^l \over (q;q)_l}\,
{q^{k(k-1)/2}\,|z|^k \over (q;q)_k}\,.
$$
Thus the double sum converges absolutely, uniformly for $z,w$ in compacta,
and it is symmetric in $z$ and $w$.\qed
\enddemo

\proclaim{Remark 2.2}\rm
Formula (2.3) is a limit case of Heine's transformation formula
$${}_2\phi_1\left[{a,b\atop c};q,z\right]
=
{(b,az;q)_\infty\over(c,z;q)_\infty}
\,
{}_2\phi_1\left[{c/b,z\atop az};q,b\right]
\tag 2.5$$
(cf\. \cite{6, (1.4.1)}).
Indeed, first let $b\to 0$, then replace $z$ by $z/a$ and let $a\to\infty$.
At least formally, by termwise limits, we obtain (2.3).
\endproclaim

\proclaim{Remark 2.3}\rm
For $w:=q^{1-n}$ ($n=1,2,\ldots$) we interpret the left-hand side of
(2.2)
by the series at its right-hand side. Then the first $n$ terms vanish, so the
summation starts with $k=n$. When we make the change of summation variable
$k=n+l$, we obtain
$$(q^{1-n};q)_\infty\,{}_1\phi_1(0;q^{1-n};q,z)=
(-1)^n\,q^{n(n-1)/2}\,z^n\,(q^{n+1};q)_\infty\,
{}_1\phi_1(0;q^{n+1};q,q^n z)
\tag 2.6$$
for $n\in\ZZ$.
(The case $n<0$ follows from the case $n>0$ of (2.6)
by changing $z$ into $q^{-n}z$.)
\endproclaim

\proclaim{Remark 2.4}\rm
Because of (2.6), the behaviour of the
two equal sides of (2.3) as $|w|\to\infty$ (cf\. (2.4))
drastically improves when $w$ runs over the values $q^{1-n}$,
$n=1,2,\ldots$. For such $w$ we
can majorize these expressions by
$$q^{n(n-1)/2}\,|z|^n\,(-|z|;q)_\infty\,(-q;q)_\infty.$$
\endproclaim

We will now restate Exton's generating function (1.5)
in terms of the notation (2.2),
and also give the short proof, for reasons of completeness.

\proclaim{Proposition 2.5}
For $z,t\in\CC$ such that $0<|t|<|z|^{-1}$ there is the absolutely
convergent expansion
$$\align
e_q(tz)\,E_q(-t^{-1}z)
=
{(t^{-1} z;q)_\infty\over(tz;q)_\infty}
&=
\sum_{n=-\infty}^\infty t^n\,z^n\,
{(q^{n+1};q)_\infty\over(q;q)_\infty}\,
{}_1\phi_1(0;q^{n+1};q,z^2)
\tag 2.7\\
&=\sum_{n=-\infty}^\infty t^n\,z^n\,
{(z^2;q)_\infty\over(q;q)_\infty}\,
{}_1\phi_1(0;z^2;q,q^{n+1}).
\tag 2.8
\endalign$$
\endproclaim
\demo{Proof}
Expansion of the left-hand side of (2.7) gives
$$\sum_{k=0}^\infty\,\sum_{l=0}^\infty
{(-1)^k\,q^{k(k-1)/2}\,t^{l-k}\,z^{l+k}\over
(q;q)_k\,(q;q)_l}
=
\sum_{k=0}^\infty\,\sum_{l=-\infty}^\infty
{(q^{l+1};q)_\infty\,(-1)^k\,q^{k(k-1)/2}\,t^{l-k}\,z^{l+k}
\over
(q;q)_\infty\,(q;q)_k}\,,
$$
which is an absolutely convergent double sum for $z,t\in\CC$
such that $0<|t|<|z|^{-1}$.
Now pass to new summation variables $k,n$ by substituting $l=k+n$.
This yields, by substitution of (2.2), the right-hand side of
(2.7).
Formula (2.8) follows by the symmetry (2.3).\qed
\enddemo

In \S3 we will show that (2.8) is a $q$-analogue of an infinite
integral of Weber and Sonine.
Replace in (2.7) $t$ by $t^{-1}$ and multiply the new identity
with the original identity. The resulting formula
$$\multline
1=\sum_{n=-\infty}^\infty\,\sum_{m=-\infty}^\infty t^{n-m}
z^n\,{(q^{n+1};q)_\infty\over(q;q)_\infty}\,
{}_1\phi_1(0;q^{n+1};q,z^2)\\
\times z^m\,{(q^{m+1};q)_\infty\over(q;q)_\infty}\,
{}_1\phi_1(0;q^{m+1};q,z^2)
\endmultline\tag 2.9$$
is absolutely convergent for $z,t\in\CC$ such that $t\ne0$ and
$|z|<|t|<|z|^{-1}$.
So equality of coefficients of equal powers of $t$ at both sides yields
a $q$-analogue of the orthogonality relations (1.8):

\proclaim{Proposition 2.6}
For $|z|<1$ and $n,m\in\ZZ$ we have
$$\multline
\sum_{k=-\infty}^\infty
z^{k+n}\, {(q^{n+k+1};q)_\infty \over (q;q)_\infty} \,
{}_1\phi_1(0;q^{n+k+1};q,z^2)\\
\times
z^{k+m}\, {(q^{m+k+1};q)_\infty \over (q;q)_\infty} \,
{}_1\phi_1(0;q^{m+k+1};q,z^2)=\de_{nm}\endmultline\tag{ 2.10}$$
and
$$\multline
\sum_{k=-\infty}^\infty
z^{k+n}\, {(z^2;q)_\infty \over (q;q)_\infty} \,
{}_1\phi_1(0;z^2;q,q^{n+k+1})\\
\times
z^{k+m}\, {(z^2;q)_\infty \over (q;q)_\infty} \,
{}_1\phi_1(0;z^2;q,q^{m+k+1})=\de_{nm},\endmultline\tag 2.11$$
where the sums at the left-hand sides are absolutely convergent,
uniformly on compact subsets of the open unit disk.
\endproclaim

Formula (2.11) follows from (2.10) by the symmetry (2.3).
In \S3 we will show that (2.11) is a $q$-version
of Hankel's Fourier-Bessel integral (1.9).

\proclaim{Remark 2.7}\rm
Analogous to Proposition 2.5 there are the two generating functions
for the $q$-Bessel functions $J_\al^{(1)}(x;q)$ and $J_\al^{(2)}(x;q)$
of Jackson and Ismail \cite{10}:
$$\multline
e_q(tz)\,e_q(-t^{-1}z)=
{1\over(tz;q)_\infty\,(-t^{-1}z;q)_\infty}\\
=\sum_{n=-\infty}^\infty t^n\,z^n\,{(q^{n+1};q)_\infty\over(q;q)_\infty}\,
{}_2\phi_1(0,0;q^{n+1};q,-z^2),\quad |z|<|t|<|z|^{-1},\endmultline\tag 2.12$$
and
$$\multline
E_q(tz)\,E_q(-t^{-1}z)=
(-tz;q)_\infty\,(t^{-1}z;q)_\infty\\
=\sum_{n=-\infty}^\infty t^n\,z^n\,
{q^{n(n-1)/2}\,(q^{n+1};q)_\infty\over(q;q)_\infty}\,
{}_0\phi_1(-;q^{n+1};q,-q^n\,z^2).\endmultline\tag 2.13$$
In a way similar to (2.10) we can now derive the
biorthogonality relations
$$\multline
\sum_{k=-\infty}^\infty
z^{k+n}\, {(q^{n+k+1};q)_\infty \over (q;q)_\infty} \,
{}_2\phi_1(0,0;q^{n+k+1};q,-z^2)\\
\times
z^{k+m}\,q^{(m+k)(m+k-1)/2}\,{(q^{m+k+1};q)_\infty \over (q;q)_\infty} \,
{}_0\phi_1(-;q^{m+k+1};q,-q^{m+k}\,z^2)=\de_{nm},
\endmultline\tag 2.14$$
valid for $|z|<1$ and $n,m\in\ZZ$.
The case $n=m$ of this result goes back to Jackson
(see also \cite{7, \S8}).
Formula (2.14) can be rewritten in several ways by
substitution of the transformations
$${}_0\phi_1(-;c;q,cz)
=(z;q)_\infty\,{}_2\phi_1(0,0;c;q,z)
={1\over(c;q)_\infty}\,{}_1\phi_1(z;0;q,c).
$$
However, this will not transform (2.14) into an orthogonality;
it remains a biorthogonality.
\endproclaim

\heading 3. Some limit transitions\endheading
Jacobi polynomials tend to Bessel functions:
$$\multline
{P_{n_N}^{(\al,\be)}(1-x^2/(2N^2))
\over
P_{n_N}^{(\al,\be)}(1)}=
{}_2F_1(-n_N,n_N+\al+\be+1;\al+1;x^2/(4N^2))\\
@> N\to\infty >>
{}_0F_1(-;\al+1;-(\la x/2)^2)
=(\la x/2)^{-\al}\,\Ga(\al+1)\,J_\al(\la x),\endmultline$$
where $n_N/N$ tends to $\la$ for $N\to\infty$.
When this limit transition is applied to the formula which recovers
a function from its Fourier-Jacobi coefficients, we obtain, at least formally,
Hankel's Fourier-Bessel integral (1.9).

The $q$-analogue of this limit transition starts with the {\sl little
$q$-Jacobi polynomials}
$$p_n(x;a,b;q):={}_2\phi_1(q^{-n},abq^{n+1}; aq;q,qx),$$
which satisfy orthogonality relations
$$\multline
{(qa,qb;q)_\infty\over(q,q^2ab;q)_\infty}\,
\sum_{k=0}^\infty(p_np_m)(q^k;a,b;q)\,(qa)^k\,
{(q^{k+1};q)_\infty\over(q^{k+1}b;q)_\infty}\\
={(qa)^n\,(1-qab)\,(qb,q;q)_n\over
(1-q^{2n+1}ab)\,(qa,qab;q)_n}\,\de_{nm},
\endmultline\tag 3.1$$
where $0<a<q^{-1}$, $b<q^{-1}$
(see Andrews and Askey \cite{2}).

It is clear that
$$p_{N-n}(q^Nx;a,b;q)=
{}_2\phi_1(q^{-N+n},abq^{N-n+1}; aq;q,q^{N+1}x)$$
tends formally (termwise) to
${}_1\phi_1(0;aq;q,q^{n+1}x)$ as $N\to\infty$.
(See Proposition A.1 for a rigorous proof of this limit result.)
Also, when we replace in (3.1) $n,m,k$ by $N-n,N-m,N+k$, respectively
(so the sum runs from $-N$ to $\infty$), and when we let $N\to\infty$,
we obtain as a formal (termwise) limit the orthogonality relations
(2.11).

In order to see that (2.11) is a $q$-analogue of Hankel's
Fourier-Bessel integral
(1.9), rewrite (2.11) as the transform pair
$$\left\{\aligned
g(q^n)&=\sum_{k=-\infty}^\infty q^{(k+n)(\al+1)}\,
{(q^{2\al+2};q^2)_\infty\over(q^2;q^2)_\infty}\,
{}_1\phi_1(0;q^{2\al+2};q^2,q^{2k+2n+2})\,f(q^k),\\
f(q^k)&=\sum_{n=-\infty}^\infty q^{(k+n)(\al+1)}\,
{(q^{2\al+2};q^2)_\infty\over(q^2;q^2)_\infty}\,
{}_1\phi_1(0;q^{2\al+2};q^2,q^{2k+2n+2})\,g(q^n),\endaligned\right.
\tag 3.2$$
where $f,g$ are $L^2$-functions on the set $\{q^k\mid k\in\ZZ\}$
with respect to counting measure.
With the ad hoc notation
$$J_\al(z;q^2):={z^\al\,(q^{2\al+2};q^2)_\infty\over (q^2;q^2)_\infty}\,
{}_1\phi_1(0;q^{2\al+2};q^2,q^2\,z^2)\tag 3.3$$
(different from Rahman's  proposal in
\cite{13, (1.13)})
and with the replacement of $f(q^k)$, $g(q^n)$ by
$q^kf(q^k)$, $q^n g(q^n)$ this becomes
$$\left\{\aligned
g(q^n)&=\sum_{k=-\infty}^\infty q^{2k}\,J_\al(q^{k+n};q^2)\,f(q^k),\\
f(q^k)&=\sum_{n=-\infty}^\infty q^{2n}\,J_\al(q^{k+n};q^2)\,g(q^n).\endaligned
\right.
\tag 3.4$$
Now observe that
$$J_\al((1-q)z;q^2)=
{(1+q)^{-\al}\,z^\al\over\Ga_{q^2}(\al+1)}\,
{}_1\phi_1(0;q^{2\al+2};q^2,(1-q^2)^2\,(qz/(1+q))^2)
$$
converges, for $q\uparrow1$, to the Bessel function
$$J_\al(z)={2^{-\al}\,z^\al\over\Ga(\al+1)}\,
{}_0F_1(-;\al+1;-z^2/4),
$$
where we used Proposition A.2 and the fact that
$${(1-q^2)^\al\,(q^{2\al+2};q^2)_\infty\over(q^2;q^2)_\infty}
={1\over\Ga_{q^2}(\al+1)}\longrightarrow{1\over\Ga(\al+1)}
$$
as $q\uparrow1$,
cf\. Andrews \cite{1, Appendix A},
Koornwinder \cite{11, Appendix B}.
We can apply this to (3.4) when we let $q\uparrow1$ under the side
condition that
${\log(1-q)\over\log q}\in2\ZZ$.
For such $q$ we can replace
$q^k$, $q^n$ in (3.4)
by $(1-q)^{1/2}$$q^k$, $(1-q)^{1/2}$$q^n$,
and next $f((1-q)^{1/2}$$q^k)$, $g((1-q)^{1/2}$$q^n)$
by $f(q^k)$, $g(q^n)$.
With the $q$-integral notation
$$\int_0^\infty h(z)\,d_qz:=(1-q)\,\sum_{j=-\infty}^\infty h(q^j)\,q^j,$$
(3.4) then takes the form
$$\left\{\aligned
g(\la)&=\int_0^\infty f( x )\,J_\al((1-q)\la x;q^2)\, x \,d_q x ,\\
f( x )&=\int_0^\infty g(\la)\,J_\al((1-q)\la x;q^2)\,\la\,d_q\la,\endaligned
\right.$$
where $\la$ in the first identity and $x$ in the second identity take the
values $q^n$, $n\in\ZZ$.
For $q\uparrow1$ we obtain, at least formally, the Hankel transform pair
$$\left\{\aligned
g(\la)&=\int_0^\infty f( x )\,J_\al(\la x)\, x \,d x ,\\
f( x )&=\int_0^\infty g(\la)\,J_\al(\la x)\,\la\,d\la,\endaligned\right.
$$
which is equivalent to (1.9).

In order to find the classical formula corresponding to (2.8),
we rewrite (2.8) first in terms of the notation (3.3):
$${(q^{\al+1-t};q^2)_\infty\over(q^{\al+1+t};q^2)_\infty}=
\sum_{n=-\infty}^\infty	q^n\,q^{nt}\,J_\al(q^n;q^2),\qquad
\Re t>-\Re\al-1.$$
For $\log(1-q)/\log q\in\ZZ$ this can be rewritten as
$${(1+q)^t\,\Ga_{q^2}((\al+1+t)/2)\over\Ga_{q^2}((\al+1-t)/2)}=
\int_0^\infty x^t\,J_\al((1-q)x;q^2)\,d_qx.$$
Formally, as $q\uparrow1$, this yields
$${2^t\,\Ga((\al+1+t)/2)\over\Ga((\al+1-t)/2)}=
\int_0^\infty x^t\,J_\al(x)\,dx,
$$
which formula, valid for $-1/2>\Re t>-\Re\al-1$,
goes back to Weber and Sonine
(cf\. \cite{15, 13.24 (1)}).

\heading 4. A $q$-analogue of Graf's addition formula\endheading
In this section we will generalize the considerations which led to
the orthogonality relations in Proposition 2.6. The resulting
formula will turn out to be a $q$-analogue of Graf's addition formula
and, at the same time, of the discontinuous integral of Weber and
Schafheitlin. In a final remark we will point out that the Graf type
addition formula part can also be done for the $q$-Bessel functions of
Jackson and Ismail, by similar methods, which go back to Heine,
the founding father of $q$-hypergeometric series.

A formula more general than (2.9) can be derived by expanding the
expression
$${(xs^{-1}t;q)_\infty\,(yt^{-1};q)_\infty\over
(yt;q)_\infty\,(xst^{-1};q)_\infty}\,.
\tag 4.1$$
as a Laurent series in $t$ ($|sx|<|t|<|y|^{-1}$) in two different ways.
On the one hand, (4.1) can be expanded by twofold substitution of the
$q$-binomial formula (cf\. \cite{6, (1.3.2)}) as:
$$\align
&{}_1\phi_0(s^{-1}xy^{-1};-;q,yt)\,{}_1\phi_0(s^{-1}yx^{-1};-;q,xst^{-1})\\
&\qquad=\sum_{k=0}^\infty\,\sum_{l=-\infty}^\infty
{(s^{-1}yx^{-1};q)_k\over(q;q)_k}\,{(s^{-1}xy^{-1};q)_l\over\,(q;q)_l}\,
s^k\,x^k\,y^l\,t^{l-k}\\
&\qquad=\sum_{n=-\infty}^\infty t^n\,y^n\,
\sum_{k=\max\{-n,0\}}^\infty
{(s^{-1}yx^{-1};q)_k\over(q;q)_k}\,
{(s^{-1}xy^{-1};q)_{n+k}\over(q;q)_{n+k}}\,(sxy)^k,\endalign$$
where we substituted $l=k+n$. Since the inner sum in the last part
can be expressed in terms of a ${}_2\phi_1$ series (depending on the sign of $n$), we obtain for $|sx|<|t|<|y|^{-1}$ the identity
$$\multline
{(xs^{-1}t,yt^{-1};q)_\infty\over
(yt,xst^{-1};q)_\infty}=
\sum_{n=0}^\infty t^n\,y^n\,
{(s^{-1}xy^{-1};q)_n\over(q;q)_n}\,
{}_2\phi_1\left[{s^{-1}yx^{-1},q^ns^{-1}xy^{-1}\atop q^{n+1}};q,sxy\right]\\
+\sum_{n=-\infty}^{-1}t^ns^{-n}x^{-n}\,
{(s^{-1}yx^{-1};q)_{-n}\over(q;q)_{-n}}\,
{}_2\phi_1\left[{s^{-1}xy^{-1},q^{-n}s^{-1}yx^{-1}\atop q^{-n+1}};q,sxy\right].
\endmultline\tag 4.2$$
We might have written (4.2) as a sum over $n$ running from $-\infty$ to $\infty$
with summand given for all integer $n$ by the left-hand side of (4.3) below,
where, for $n<0$,
an interpretation of the ${}_2\phi_1$ similar to our
convention in Remark 2.3 is used. Thus the analogue of (2.6) becomes:
$$\multline
t^ny^n\,{(s^{-1}xy^{-1},q^{1+n};q)_\infty\over(q^ns^{-1}xy^{-1},q;q)_\infty}\,
{}_2\phi_1\left[{s^{-1}yx^{-1},q^ns^{-1}xy^{-1}\atop q^{1+n}};q,sxy\right]\\
=t^ns^{-n}x^{-n}\,
{(s^{-1}yx^{-1},q^{1-n};q)_\infty\over(q^{-n}s^{-1}yx^{-1},q;q)_\infty}\,
{}_2\phi_1\left[{s^{-1}xy^{-1},q^{-n}s^{-1}yx^{-1}\atop q^{1-n}};q,sxy\right].
\endmultline\tag 4.3$$
Here we should exclude the case that both $1+n$ and
$n+{}^q\log(s^{-1}xy^{-1})$ are non-positive integers, or that
both $1-n$ and $-n+{}^q\log(s^{-1}yx^{-1})$ are non-positive integers.

Note that (4.2) reduces to (2.7) in the
special case $x=0$ (and also, in view of (2.6), for $y=0$).
Formula (4.2) is a $q$-analogue of
(1.4) with $z$ and $t$ replaced by 
$2(y-s^{-1}x)^{1/2}\*(y-sx)^{1/2}$ and $t(y-s^{-1}x)^{1/2}\*(y-sx)^{-1/2}$,
respectively.

On the other hand we expand (4.1) by twofold substitution of
(2.7):
$$\multline
{(xs^{-1}t;q)_\infty\,(yt^{-1};q)_\infty\over
(yt;q)_\infty\,(xst^{-1};q)_\infty}
=\sum_{n=-\infty}^\infty\,\sum_{k=-\infty}^\infty t^n\,s^k\\
\times y^{n+k}\,
{(q^{n+k+1};q)_\infty\over(q;q)_\infty}\,{}_1\phi_1(0;q^{n+k+1};q,y^2)\,
x^k\,{(q^{k+1};q)_\infty\over(q;q)_\infty}\,{}_1\phi_1(0;q^{k+1};q,x^2),
\endmultline\tag 4.4$$
which generalizes (2.9).
When we compare coefficients of equal powers of $t$ in (4.2) and
(4.4), we obtain:

\proclaim{Proposition 4.1}
For $|sxy|<1$ we have
$$\multline
\sum_{k=-\infty}^\infty s^k\,y^{n+k}\,
{(q^{n+k+1};q)_\infty\over(q;q)_\infty}\,{}_1\phi_1(0;q^{n+k+1};q,y^2)\,
x^k\,{(q^{k+1};q)_\infty\over(q;q)_\infty}\,{}_1\phi_1(0;q^{k+1};q,x^2)\\
=\left\{\aligned
&y^n\,
{(s^{-1}xy^{-1};q)_n\over(q;q)_n}\,
{}_2\phi_1\left[{s^{-1}yx^{-1},q^ns^{-1}xy^{-1}\atop q^{n+1}};q,sxy\right]\quad
\text{\rm if $n\ge0$,}\\
&s^{-n}x^{-n}\,
{(s^{-1}yx^{-1};q)_{-n}\over(q;q)_{-n}}\,
{}_2\phi_1\left[{s^{-1}xy^{-1},q^{-n}s^{-1}yx^{-1}\atop q^{-n+1}};q,sxy\right]\quad
\text{\rm if $n<0$.}\endaligned\right.
\endmultline\tag 4.5$$
\endproclaim

Formula (4.5) is a $q$-analogue of the addition formula
$$\sum_{k=-\infty}^\infty s^k\,J_{n+k}(y)\,J_k(x)=
\left({y-s^{-1}x\over y-sx}\right)^{n/2}\,
J_n\left(\sqrt{(y-s^{-1}x)(y-sx)}\right),
\tag 4.6$$
due to Graf (cf\. \cite{15, \S11.3 (1)}).
The special case $n=0$ of (4.5) is a $q$-analogue of Neumann's
addition formula for Bessel functions $J_0$
(cf\. \cite{15, \S11.2 (1)}).
In the special case $x=y$, $s=1$ the right-hand side of (4.5) becomes
$y^n\,\de_{n,0}$, so then (4.5) reduces to the orthogonality
relations (2.10).

When we apply the symmetries (2.5) and (4.3) to the left-hand
side of (4.5) and also apply (2.3) to the right-hand side of
(4.5), and replace $n$ by $n-m$ and next $k$ by $k+m$ then
we obtain an equivalent identity:

\proclaim{Proposition 4.2}
For $|sxy|<1$ and $x^2,y^2\ne q^{-l}$ {\rm($l=0,1,2,\ldots$)},
we have
$$\multline
\sum_{k=-\infty}^\infty s^k\,y^{n+k}\,{(y^2;q)_\infty\over(q;q)_\infty}\,
{}_1\phi_1(0;y^2;q,q^{n+k+1})\,
x^{m+k}\,{(x^2;q)_\infty\over(q;q)_\infty}\,{}_1\phi_1(0;x^2;q,q^{m+k+1})\\
=\left\{\aligned
&s^{-m}\,y^{n-m}\,{(s^{-1}xy^{-1},y^2;q)_\infty\over(sxy,q;q)_\infty}\,
{}_2\phi_1\left[{qsx^{-1}y,sxy\atop y^2};q,q^{n-m}s^{-1}xy^{-1}\right],\\
&s^{-n}\,x^{m-n}\,{(s^{-1}yx^{-1},x^2;q)_\infty\over(sxy,q;q)_\infty}\,
{}_2\phi_1\left[{qsxy^{-1},sxy\atop x^2};q,q^{m-n}s^{-1}yx^{-1}\right],\endaligned\right.
\endmultline\tag 4.7$$
where the first expression on the right-hand side
can be chosen if
$s^{-1}xy^{-1}$ and $q^{n-m}s^{-1}xy^{-1}$ are not equal
to a non-positive integer power of $q$, and where the second
expression on the right-hand side can be chosen if
$s^{-1}yx^{-1}$ and $q^{m-n}s^{-1}yx^{-1}$ are not equal to a non-positive integer
power of $q$.
Thus in the generic case both expressions on the right-hand side are equal.
\endproclaim

See Appendix B for the detailed proof.
Observe the special case $x=y$ of (4.7):
$$\multline
\sum_{k=-\infty}^\infty s^k\,z^{n+k}\,{(z^2;q)_\infty\over(q;q)_\infty}\,
{}_1\phi_1(0;z^2;q,q^{n+k+1})\,
z^{m+k}\,{(z^2;q)_\infty\over(q;q)_\infty}\,{}_1\phi_1(0;z^2;q,q^{m+k+1})\\
=\left\{\aligned
&s^{-m}\,z^{n-m}\,{(s^{-1},z^2;q)_\infty\over(sz^2,q;q)_\infty}\,
{}_2\phi_1\left[{qs,sz^2\atop z^2};q,q^{n-m}s^{-1}\right],\\
&s^{-n}\,z^{m-n}\,{(s^{-1},z^2;q)_\infty\over(sz^2,q;q)_\infty}\,
{}_2\phi_1\left[{qs,sz^2\atop z^2};q,q^{m-n}s^{-1}\right],\endaligned\right.
\endmultline\tag 4.8$$
where $|sz^2|<1$, $z^2$ is not a non-positive integer power of $q$
and the expression on the right-hand side has to be chosen similarly as in
Proposition 4.2. This can be considered as a kind of
Poisson kernel for the orthogonal
system with orthogonality relations (2.11).
Inspection of the right-hand side of (4.8) shows that this
kernel is
positive if $0<z<1$ and $1<s<\min\{q^{-1},z^{-2}\}$.

Let us look for the classical analogue of formula (4.7).
In (4.7) first replace $q$ by $q^2$, then
$x,y,s$ by $q^{\al+1},q^{\be+1},q^{-\ga-1}$, respectively.
Then, with the notation (3.3), formula
(4.7) can be rewritten as
$$\multline
\sum_{k=-\infty}^\infty q^{k(-\ga+1)}\,J_\al(q^{m+k};q^2)\,J_\be(q^{n+k};q^2)
\\
=\left\{\aligned
q^{n\be}\,q^{m(\ga-\be-1)}\,&
{(q^{\al-\be+\ga+1},q^{2\be+2};q^2)_\infty\over
(q^{\al+\be-\ga+1},q^2;q^2)_\infty}\\
&\times{}_2\phi_1(q^{\be-\al-\ga+1},q^{\al+\be-\ga+1}; q^{2\be+2};
q^2,q^{2n-2m+\al-\be+\ga+1}),\\
q^{m\al}\,q^{n(\ga-\al-1)}\,&
{(q^{\be-\al+\ga+1},q^{2\al+2};q^2)_\infty\over
(q^{\al+\be-\ga+1},q^2;q^2)_\infty}\\
&\times{}_2\phi_1(q^{\al-\be-\ga+1},q^{\al+\be-\ga+1}; q^{2\al+2};
q^2,q^{2m-2n+\be-\al+\ga+1}),\endaligned\right.\endmultline$$
where $\Re(\al+\be-\ga+1)>0$,
$\al,\be$ are not equal to a negative integer, and where we can choose
the first expression on the right-hand side if $(\alpha-\beta+\gamma+1)/2$ and
$n-m+(\alpha-\beta+\gamma+1)/2$ are not a non-positive integer
and we can choose
the second expression on the right-hand side if $(\beta-\alpha+\gamma+1)/2$ and
$m-n+(\beta-\alpha+\gamma+1)/2$ are not a non-positive integer.

Now replace $q^k$ by $q^k(1-q)$ (with
${\log(1-q)\over\log q}\in\ZZ$) and let $m,n$ depend on $q$ such that,
as $q\uparrow1$, $q^m$ tends to $a$ and $q^n$ tends to $b$.
Depending on whether $b<a$ or $b>a$, make the formal limit transition
$q\uparrow1$ in the first or second identity, respectively.
Then, for $\Re(\al+\be-\ga+1)>0$, $\Re\ga>-1$,
we obtain the
discontinuous integral of Weber and Schafheitlin:
$$\multline
2^\ga\,\int_0^\infty x^{-\ga}\,J_\al(ax)\,J_\be(bx)\,dx\\
=\left\{\aligned
{a^{\ga-\be-1}\,b^\be\,
\Ga((\al+\be-\ga+1)/2)\over
\Ga((\al-\be+\ga+1)/2)\,\Ga(\be+1)}&\\
\times{}_2F_1((\be-\al-\ga+1)/2,&(\al+\be-\ga+1)/2;\be+1;
b^2/ a^2)\qquad
\text{if $ b<a$,}\\
{a^\al\,b^{\ga-\al-1}\,
\Ga((\al+\be-\ga+1)/2)\over
\Ga((\be-\al+\ga+1)/2)\,\Ga(\al+1)}&\\
\times{}_2F_1((\al-\be-\ga+1)/2,&(\al+\be-\ga+1)/2;\al+1;
a^2/b^2)\qquad
\text{if $a<b$,}\endaligned\right.\endmultline
$$
(cf\. \cite{15, \S13.4 (2)}).
Note that the two analytic expressions at the right-hand side are no longer
equal, as they were in the $q$-case.

\proclaim{Remark 4.3}\rm
Analogous to (4.2) we have
$$\multline
{(xs^{-1}t;q)_\infty\,(-xst^{-1};q)_\infty
\over
(yt;q)_\infty\,(-yt^{-1};q)_\infty}=
\sum_{n=-\infty}^\infty
t^n\,y^n\,{(s^{-1}xy^{-1},q^{n+1};q)_\infty\over
(q^ns^{-1}xy^{-1},q;q)_\infty}\\
\times
{}_2\phi_1\left[{q^ns^{-1}xy^{-1},sxy^{-1}\atop q^{n+1}};q,-y^2\right],
\quad|y|<|t|<|y|^{-1}.\endmultline
\tag 4.9$$
The case $s=1$ of (4.9) goes
back to Heine \cite{9, p.121} (see also
Hahn \cite{7, \S8}).
Its special cases $(x,y,s)=(0,z,1)$ and $(-z,0,1)$ are the formulas
(2.12) and (2.13).
Like (4.2), formula (4.9) is a $q$-analogue of
(1.4) with $z$ and $t$ replaced by 
$2(y-s^{-1}x)^{1/2}\*(y-sx)^{1/2}$ and $t(y-s^{-1}x)^{1/2}\*(y-sx)^{-{1/2}}$,
respectively.
Similarly to (4.5) we obtain from (2.12), (2.13) and
(4.9) that, for $|y|<1$:
$$\multline
y^n\,{(s^{-1}xy^{-1},q^{n+1};q)_\infty\over(q^ns^{-1}xy^{-1},q;q)_\infty}\,
{}_2\phi_1\left[{q^ns^{-1}xy^{-1},sxy^{-1}\atop q^{n+1}};q,-y^2\right]
=
\sum_{k=-\infty}^\infty s^k\,y^{k+n}\\
\times{(q^{k+n+1};q)_\infty\over(q;q)_\infty}\,
{}_2\phi_1(0,0;q^{k+n+1};q,-y^2)\,
x^k\,q^{k(k-1)/2}\,{(q^{k+1};q)_\infty\over(q;q)_\infty}\,
{}_0\phi_1(-;q^{k+1};q,-q^k x^2).\endmultline
\tag 4.10$$
The special case $n=0$ of this formula is the limit case
$\nu\downarrow0$ of Rahman's addition
formula \cite{12, (1.10)}.
Like (4.5), formula (4.10) is a $q$-analogue of Graf's
addition formula (4.6).
The special case $x=y$, $s=1$ of (4.10) gives the
biorthogonality relations (2.14).
\endproclaim

\heading 5. $q$-Analogues of the Fourier-cosine and Fourier-sine
transform\endheading
Put
$$\align
\cos(z;q^2):=&{}_1\phi_1(0;q;q^2,q^2\,z^2)\\
=&\sum_{k=0}^\infty{(-1)^k\,q^{k(k+1)}\,z^{2k}\over(q;q)_{2k}}
\tag 5.1\\
=&{(q^2;q^2)_\infty\over(q;q^2)_\infty}\,z^{1/2}\,J_{-1/2}(z;q^2)
\tag 5.2\endalign$$
and
$$\align
\sin(z;q^2):=&(1-q)^{-1}\,z\,{}_1\phi_1(0;q^3;q^2,q^2\,z^2)\\
=&\sum_{k=0}^\infty{(-1)^k\,q^{k(k+1)}\,z^{2k+1}\over(q;q)_{2k+1}}
\tag 5.3\\
=&{(q^2;q^2)_\infty\over(q;q^2)_\infty}\,z^{1/2}\,J_{1/2}(z;q^2).
\tag 5.4\endalign$$
Here we have used the notation (3.3).
(The functions introduced above should not be confused with the functions
$\cos_q$ and $\sin_q$ considered in \cite{6, Exercise 1.14}.)
Clearly we have the formal (termwise) limits
$$\cos((1-q)z;q^2)\longrightarrow\cos z\quad\hbox{and}\quad
\sin((1-q)z;q^2)\longrightarrow\sin z$$
as $q\uparrow1$. By Proposition A.2 these limit transitions hold pointwise,
uniformly on compacta.
When we substitute (5.2) or (5.4) in (3.4) and replace
$f(q^k)$, $g(q^n)$ by $q^{-k/2}$$f(q^k)$, $q^{-n/2}$$g(q^n)$,
we obtain the transform pairs
$$\left\{\aligned
g(q^n)&={(q;q^2)_\infty\over(q^2;q^2)_\infty}\,
\sum_{k=-\infty}^\infty q^k\,
\left\{\matrix\cos(q^{k+n};q^2)\\\text{or}\\\sin(q^{k+n};q^2)\endmatrix
\right\}\,
f(q^k),\\
f(q^k)&={(q;q^2)_\infty\over(q^2;q^2)_\infty}\,
\sum_{n=-\infty}^\infty q^n\,
\left\{\matrix\cos(q^{k+n};q^2)\\\text{or}\\\sin(q^{k+n};q^2)\endmatrix
\right\}\,
g(q^n).\endaligned\right.
\tag 5.5$$

The transformations $f\mapsto g$ and $g\mapsto f$ of (5.5)
establish an isometry of Hilbert spaces:
$$\sum_{k=-\infty}^\infty q^k\,|f(q^k)|^2
= \sum_{n=-\infty}^\infty q^n\,|g(q^n)|^2.
$$

Now let $q\uparrow1$ under the side
condition that
${\log(1-q)\over\log q}\in2\ZZ$.
Replace
$q^k$, $q^n$ in (5.5)
by $(1-q)^{1/2}\*q^k$, $(1-q)^{1/2}\*q^n$,
and then $f((1-q)^{1/2}\*q^k)$, $g((1-q)^{1/2}\*q^n)$
by $f(q^k)$, $g(q^n)$.
Then (5.5) takes the form
$$\left\{\aligned
g(\la)&={(1+q)^{1/2}\over\Ga_{q^2}(1/2)}\,
\int_0^\infty f(x)\,
\left\{\matrix\cos((1-q)\la x;q^2)\\\text{or}\\\sin((1-q)\la x;q^2)\endmatrix
\right\}
\,d_qx,\\
f(x)&={(1+q)^{1/2}\over\Ga_{q^2}(1/2)}\,
\int_0^\infty g(\la)\,
\left\{\matrix\cos((1-q)\la x;q^2)\\\text{or}\\\sin((1-q)\la x;q^2)\endmatrix
\right\}
\,d_q\la.\endaligned\right.
$$
Formally, as $q\uparrow 1$, we obtain the classical Fourier pairs
$$\eqalignno{
g(\la)=\sqrt{2/\pi}\,\int_0^\infty f(x)\,\cos(\la x)\,dx,\qquad&
f(x)=\sqrt{2/\pi}\,\int_0^\infty g(\la)\,\cos(\la x)\,d\la\cr
\noalign{\hbox{and}}
g(\la)=\sqrt{2/\pi}\,\int_0^\infty f(x)\,\sin(\la x)\,dx,\qquad&
f(x)=\sqrt{2/\pi}\,\int_0^\infty g(\la)\,\sin(\la x)\,d\la.\cr}
$$

With the notation
$$(D_qf)(z)=D_{q,z}\,f(z):=
{f(z)-f(qz)\over(1-q)z}
$$
for the $q$-derivative, we obtain from (5.1) and (5.3) that
$$\eqalign{
(1-q)\,D_{q,z}\,\cos(z;q^2)&=-q\,\sin(qz;q^2),\cr
(1-q)\,D_{q,z}\,\sin(z;q^2)&=    \cos( z;q^2).\cr}
$$
Hence
$$(1-q)^2\,(D_q^2f)(q^{-1}z)=\cases
-q^2\,\la^2\,f(z)&\text{if $f(z)=\cos(\la z;q^2)$,}\\
-q  \,\la^2\,f(z)&\text{if $f(z)=\sin(\la z;q^2)$.}\endcases
$$
So the two systems of functions $z\mapsto\cos(q^nz)$ ($n\in\ZZ$)
and $z\mapsto\sin(q^nz)$ ($n\in\ZZ$) have disjoint eigenvalues with
respect to the operator which sends $f$ to the function
$$z\mapsto(1-q)^2\,(D_q^2f)(q^{-1}z)=
q\,z^{-2}\bigl(q\,f(q^{-1}z)-(1+q)\,f(z)+f(qz)\bigr).
$$
This operator has also the selfadjointness property
$$\sum_{k=-\infty}^\infty q^k\,(D_q^2f)(q^{k-1})\,g(q^k)=
\sum_{k=-\infty}^\infty q^k\,f(q^k)\,(D_q^2g)(q^{k-1})
$$
for $f$, $g$ of finite support on $\{q^k\mid k\in\ZZ\}$.

Observe that the $q$-deformation of $d^2/dx^2$ considered above
yields a symmetry breaking. The two-dimensional eigenspaces of
$d^2/dx^2$ are broken apart into one-dimensional eigenspaces.
Therefore it does not seem to be very useful to consider a $q$-exponential
built from the functions defined by (5.1) and (5.3).
Any linear combination $f(z)$ of $\cos(\la z;q^2)$ and $\sin(\la z;q^2)$
will no longer satisfy an eigenfunction equation
$$(1-q)^2\,(D_q^2f)(q^{-1}z)=\mu\,f(z),\tag 5.6$$
while the nice function
$$f(z):={}_1\phi_1(0;-q^{1/2};q^{1/2},\pm iq^{3/4}\la z)=
\cos(\la z;q^2)\mp iq^{1/4}\sin(q^{1/2}\la z;q^2)$$
(cf\. Exton \cite{4, 5.2.2.1}),
which satisfies (5.6) with $\mu=-q^2\la^2$, no longer remains within
the spectral decompositions implied by (5.5).

\heading Appendix A. Rigorous proofs of some limit results\endheading
In this Appendix we will give proofs of the limit transitions from
little $q$-Jacobi polynomials to $q$-Bessel functions and from $q$-Bessel
functions to ordinary Bessel functions.
\proclaim{Proposition A.1}
For $0<a<q^{-1}$, $0\le b<q^{-1}$ we have
$$\lim_{n\to\infty}
{}_2\phi_1(q^{-n},q^{n+1}ab; qa;q,q^n x)
=
{}_1\phi_1(0; qa;q,x),$$
uniformly for $x$ in compact subsets of $\CC$.
\endproclaim
\demo{Proof}
Put
$$\multline
R_n(x):=
{}_2\phi_1(q^{-n},q^{n+1}ab; qa;q,q^n x)
-
{}_1\phi_1(0; qa;q,x)\\
=
\sum_{k=1}^\infty{(-1)^k\,q^{ k(k-1)/2}\,x^k\over
(qa;q)_k\,(q;q)_k}\,
\biggl(-1+\prod_{j=1}^k(1-q^{n-j+1})\,(1-q^{n+j}ab)\biggr).\endmultline$$
Since
$$\prod_{j=1}^k(1-x_j)\ge1-\sum_{j=1}^k x_j\qquad
\text{if\quad$0\le x_j\le1$, $j=1,\dots,k$,}$$
we have
$$\aligned
\biggl|-1+\prod_{j=1}^k(1-q^{n-j+1})\,(1-q^{n+j}ab)\biggr|
\le&\sum_{j=1}^k q^{n-j+1}+\sum_{j=1}^k q^{n+j}ab\\
\le&{q^{n-k+1}(1+ab)\over 1-q}\endaligned$$
if $k\le n$, while the left-hand side is equal to $1$ otherwise.
Thus, for $|x|\le M$:
$$|R_n(x)|\le{q^{n+1}(1+ab)\over1-q}\,
\sum_{k=1}^\infty{q^{ k(k-1)/2}\,q^{-k}M^k\over
(qa;q)_k\,(q;q)_k}+\sum_{k=n+1}^\infty{q^{ k(k-1)/2}\,M^k\over
(qa;q)_k\,(q;q)_k}\,,
$$
which tends to 0 as $n\to\infty$.\qed
\enddemo

\proclaim{Proposition A.2}
For $\al>-1$ we have
$$\lim_{q\uparrow1}\,{}_1\phi_1(0;q^{\al+1};q,(1-q)^2z)=
{}_0F_1(-;\al+1;-z),$$
uniformly for $z$ in compact subsets of $\CC$.
\endproclaim
\demo{Proof}
$${}_1\phi_1(0;q^{\al+1};q,(1-q)^2z)=
\sum_{k=0}^\infty{(-1)^k\,q^{ k(k-1)/2}\,(1-q)^{2k}\,z^k\over
(q^{\al+1};q)_k\,(q;q)_k}$$
and the summand in the sum at the right-hand side can be majorized by
$$(q^{-\al/2}|z|)^k\,\prod_{j=0}^{k-1}
{(q^{(\al+j)/2}-q^{1+(\al+j)/2})\,(q^{j/2}-q^{1+j/2})\over
(1-q^{1+\al+j})\,(1-q^{1+j})}\,.
\tag A.1$$
Now, by \cite{11, Lemma A.1}
(read $-1\le\mu-\la$ instead of $0\le\mu-\la$ in the formulation of that
lemma), we see that
$${q^{(\al+j)/2}-q^{1+(\al+j)/2}\over1-q^{1+\al+j}}$$
increases to $(1+\al+j)^{-1}$ as $q\uparrow1$
if $\al+j\ge0$.
So, the expression in (A.1), for $1/2<q<1$, is dominated by
$$\eqalignno{
{(2^{\al/2}|z|)^k\over(\al+1)_k\,k!}\quad&
\hbox{if $\al\ge0$}\cr
\noalign{\hbox{and by}}
\const{|z|^k\over(\al+1)_k\,k!}\quad&
\hbox{if $-1<\al<0$.}\cr}
$$
So the proposition follows by dominated convergence.\qed
\enddemo

\heading Appendix B. Proof of Propostion 4.2\endheading
First observe that
$u(z):={}_2\phi_1(a,b;c;q,z)$ ($c$ not a non-positive integer power of
$q$) is analytic in $z$ on $\CC$ outside the set $\{q^{-l}\mid l=0,1,2,\ldots\}$,
on which points the function has possible poles.
This follows immediately by rewriting
the second order $q$-difference equation for $u$, given in 
[6, Exercise 1.13], as
$$
u(z)=(z-1)^{-1}\big(((q^a+q^b)z-q^{c-1}-1)u(qz)+(q^{c-1}-q^{a+b}z)u(q^2z)\big).
\tag B.1
$$
Indeed, $u(z)$ is analytic in $z$ for $|z|<1$. Then (B.1) shows that $u(z)$ can
be extended to an analytic function in $z$ for $|z|<q^{-1}$
with possible pole in 1,
and then by iteration to an analytic function for $|z|<q^{-l}$ with possible poles
in $1,q^{-1},\ldots,q^{-l+1}$.

Now use the identity (see [6, (4.3.2)])
$$\multline
{}_2\phi_1(a,b;c;q,z)
={(b,c/a,az,q/(az);q)_\infty\over (c,b/a,z,q/z;q)_\infty}\,
{}_2\phi_1(a,aq/c;aq/b;q,cq/(abz))\\
+{(a,c/b,bz,q/(bz);q)_\infty\over(c,a/b,z,q/z;q)_\infty}\,
{}_2\phi_1(b,bq/c;bq/a;q,cq/(abz)),
\endmultline \tag B.2$$
where $a/b$ and $z$ are not an integer power of $q$, where $c$ and $cq/(abz)$
are not a non-positive integer power of $q$, and where $a,b,c,z\ne0$.
Now put $z:=q^{n+1}/a$ ($n=0,1,2,\ldots$) in (B.2). Then the first term
on the right-hand side of
(B.2) vanishes and after some rewriting of the coefficient in the second term
we obtain
$$
{}_2\phi_1(a,b;c;q,q^{n+1}/a)=
{(qb/a,c/b;q)_\infty\over(c,q/a;q)_\infty}\,b^{-n}\,
{}_2\phi_1(qb/c,b;qb/a;q,q^{-n}c/b),
\tag B.3$$
where initially
$a/b$ and $a$ are not an integer power of $q$, $c$ and $q^{-n}c/b$
are not a non-positive integer power of $q$, and $a,b,c\ne0$.
By continuity we see that (B.3) more generally holds if
$c,qb/a,q/a,q^{-n}c/b$ are not a non-positive integer power of $q$ and $a,b,c\ne0$.
The case of (B.3) with $n$ a non-positive integer then follows by redefining the
parameters in (B.3). Finally we can state:
\proclaim{Lemma B.1}
Formula {\rm (B.3)} is valid for integer $n$ if
$c,qb/a,q/a,q^{n+1}/a,q^{-n}c/b,c/b$ are not a non-positive integer power of $q$
and $a,b,c\ne0$.
\endproclaim

Now consider (2.5). Initially this is proved (see [6, (1.4.1)]) for $|z|,|b|<1$ and
$c,az$ not a non-positive integer power of $q$. By analytic continuation
(2.5) is more generally valid for $z,b,c,az$ not a non-positive integer power of $q$.

Let us turn to the proof of Proposition 4.2. Proceed with (4.5) as
written just before Proposition 4.2 and take into account the constraints for
(2.5) as just mentioned. Then we initially obtain, under the assumptions that
$|sxy|<1$ and $x^2,y^2$ are not non-positive integer powers of $q$, that
(4.7) holds with the first expression taken on the right-hand side if
$n\ge m$ and $q^{n-m}s^{-1}xy^{-1}$ not a non-positive integer power of $q$,
and that it holds with the second expression taken on the right-hand side if
$m\ge n$ and $q^{m-n}s^{-1}yx^{-1}$ not a non-positive integer power of $q$.
However, by Lemma B.1, the two expressions on the right-hand side of (4.7)
are equal to each other for all integer values of $m,n$ if
$s^{-1}xy^{-1}, q^{n-m}s^{-1}xy^{-1},s^{-1}yx^{-1},q^{m-n}s^{-1}yx^{-1}$ are not a
non-positive integer power of $q$.
Hence, under these last constraints the left-hand side of (4.7) is equal
to both expressions on the right-hand side for all integer values of $m,n$.
Finally, if $s^{-1}xy^{-1}, q^{n-m}s^{-1}xy^{-1}$ are not a non-positive integer power
of $q$ but $s^{-1}yx^{-1}$ or $q^{m-n}s^{-1}yx^{-1}$ is a non-positive integer power
of $q$, then by continuity the left-hand side of (4.7) is still equal to the
first expression on the right-hand side. Similarly, if
$s^{-1}yx^{-1},q^{m-n}s^{-1}yx^{-1}$  are not a non-positive integer power
of $q$ but $s^{-1}xy^{-1}$ or $ q^{n-m}s^{-1}xy^{-1}$ is a non-positive integer power
of $q$, then the left-hand side of (4.7) is still equal to the
second expression on the right-hand side.

\enddocument